\documentclass{amsart}
\newcommand{\note}{\noindent {\bf Notation. }}

\newcommand{\ws}{\hspace{4pt}}
\newtheorem{theorem}{Theorem}

\newtheorem{proposition}{Proposition}

\newtheorem{lemma}{Lemma}

\newtheorem{defi}{Definition}
\usepackage[]{amssymb,amsopn}
\usepackage[]{amsmath}
\usepackage[]{eucal}
\usepackage[]{latexsym}
\usepackage{color}

\makeatletter
\@namedef{subjclassname@2020}{%
  \textup{2020} Mathematics Subject Classification}
\makeatother

\begin{document}

\title[Removable Sets]{Removable Sets for  Fractional Heat and Fractional Bessel-Heat Equations}
\author{Mouna Chegaar and \'A. P. Horv\'ath}

\subjclass[2020]{35R11, 31C15}
\keywords{Fractional heat equations, $L^p$-removable sets, capacity, Sobolev-type spaces, Hankel transform}

\begin{abstract}
We examine the fractional heat diffusion equations $L_{\gamma,a}:=(-\Delta_a)^{\frac{\gamma}{2}}+\partial_t$, where $\Delta_a$ is the Laplace- or the Bessel-Laplace operator. We give conditions  for removability which are sufficient and which are necessary, by $L^p$-capacities. Introducing a spherical modulus of smoothness we can treat the Laplace and Bessel-Laplace cases together.
\end{abstract}
\maketitle

\section{Introduction}

Removable sets, with respect to a partial differential operator, are sets to which the solution can be extended in a certain sense. Characterization of removability of sets with respect to harmonic functions by its Hausdorff dimension is due to Carleson, see \cite{c}. After this result, the study was extended to several types of partial differential equations, and the criteria for removability were determined using capacity or Hausdorff measure. As an illustration we mention the following classical result, cf. \cite[Theorem 2.7.4 and Corollary 3.3.4]{ah}.

\medskip

\noindent {\bf Theorem A} {\it Let $\mathcal{L}$ be an elliptic linear partial differential operator of order $ \alpha< n$ with constant coefficients, and let $K\subset \mathbb{R}^n$ be compact. Then $K$ is removable for $\mathcal{L}$ in $L^p$, $1<p<\infty$, if and only if $C_{\alpha,p'}(K)=0$.}

\medskip

\noindent $C_{\alpha,p'}$ is the $L^{p'}$-capacity, which is given by Sobolev space norm. A similar result can be proven for "Bessel-elliptic" equations with Bessel $p$-capacity, see \cite[Theorem 1]{h}. Recently, conditions for removability have been given for certain fractional heat-diffusion equations in Lipschitz-class case, see \cite{mpt}, in $L^\infty$ case, see \cite{mp} and for Cantor sets, see \cite{he}. In the last two papers the main results are proved in $\frac{1}{2}$-fractional case. Our purpose is to get similar results to exponents from $(0,1)$ and in $L^p$ with $1<p<\infty$.

\medskip

\noindent The second motivation for this paper is to emphasize the close connection between the Laplace and Bessel-Laplace operators, where Bessel-Laplace operator is $$\Delta_a=\sum_{i=1}^n(\partial_{x_ix_i}+\frac{a_i}{x_i}\partial_{x_i}), \ws \ws a_i>0.$$
One way to describe this relationship is through the Poisson operator, see \cite[Statement 4, page 137]{ss}.
Our approach is different.
Technically, the Fourier transform with standard convolution in  Laplace case behaves like the Hankel transform with Bessel convolution in Bessel-Laplace case. This coincidence allows the homogeneous discussion. \\
The other important ingredient is the observation that taking average on the sphere reduces the problems to one variable ones, cf. Definition \ref{sd}. The method of reducing problems to radial ones is useful in applications. One of the most important tools in this area is Bessel convolution and Hankel transform, see e.g. \cite{m}, \cite{hm}.\\ Below we homogenize the two problems with the defined "spherical difference" via a one-dimensional  "modulus of smoothness  and K-functional" type  chain of ideas.
This leads to the uniformly formulated result of Proposition \ref{lp}.\\ Of course, there are some differences between the two cases. The first one is that the behavior of the Riesz kernel is different from the Bessel one at infinity. Thus the formulation of Theorem \ref{th1} refers to the Bessel case.\\
On the other hand, Bessel translation has some inconvenient properties. For instance it does not commute with the first partials and with product. Taking into consideration these differences, in Laplace case the statement of Proposition \ref{lp} can be extended to all exponents, cf. Theorem \ref{t2}.

\section{Notation, preliminaries}

In this preparatory section we introduce the notation and collect the results which we will use subsequently. Since we will treat the Laplace and Bessel cases together, we will give a more or less homogeneous notation.

\subsection{Spaces} In Bessel case, we introduce the multiindex $a:=a_1,\dots ,a_n$ where $a_i=2\alpha_i+1$ and $\alpha_i>-\frac{1}{2}$, $i=1, \dots , n$. $|a|=\sum_{i=1}^na_i$. In Laplace case $|a|=0$.\\
$\mathbb{R}^n_{(+)}$ means $\mathbb{R}^n_+$, i.e. the (open) positive orthant in Bessel case, and $\mathbb{R}^n$ in Laplace case. $(x,t)\in \mathbb{R}^n_{(+)}\times \mathbb{R}_+$. On the positive orthant we need the positive parts of sets: $H_+:=H\cap \mathbb{R}^n_+$. \\
The weighted space $L^p_a(\mathbb{R}^n_{(+)})$ is defined by the norm
$$\|f\|_{p,a}^p=\int_{\mathbb{R}^n_{(+)}}|f(x)|^px^adx,$$
which is just the standard $L^p$ space if $|a|=0$.\\
On $\mathbb{R}^{n+1}_{(+)}$ we introduce the $L^p_A(\mathbb{R}^{n+1}_{(+)})$ space with norm
$$\|f\|_{p,A}^p=\int_0^\infty\int_{\mathbb{R}^n_{(+)}}|f(x)|^px^adxdt.$$
If $|a|=0$, $\mathcal{S}_{(e)}=\mathcal{S}(\mathbb{R}^n)$ is the Schwartz class.
If $|a|>0$, $\mathcal{S}_{(e)}$ stands for the even functions from the Schwartz class, i.e.\\ $\mathcal{S}_{(e)}:=\left\{f\in \mathcal{S} : \frac{\partial^{2m+1}}{\partial x_i^{2m+1}}f\large{|}_{x_i=0}=0, \ws \forall m \in \mathbb{N}, \ws i=1, \dots , n\right\}.$  $\mathcal{S}_t$ denotes the Schwartz functions on $\mathbb{R}^{n+1}$ restricted to the $t>0$ half-space. $\mathcal{S}_{(e),t}$ contains that functions from $\mathcal{S}_{t}$ which are in $\mathcal{S}_{(e)}$ in the first $n$ variables.\\
$\mathcal{S}_{(e)}'$, $\mathcal{S}_{(e),t}'$ ... stand for the corresponding distributions.\\
The scalar products referring the introduced spaces are
$$\langle u,\varphi\rangle_a=\int_{\mathbb{R}^n_{(+)}}u(x)\varphi(x)x^adx \ws\ws \mbox{and} \ws \ws \langle u,\varphi\rangle_A=\int_0^\infty \int_{\mathbb{R}^n_{(+)}}u(x,t)\varphi(x,t)x^adxdt.$$
In distributional sense we use the same notation. $\delta_a$ $(\delta_A)$ stands for the weighted Dirac delta distribution, for which $\langle \delta_a,\varphi\rangle_a=\varphi(0)$ and $\langle \delta_A,\varphi\rangle_A=\varphi(0,0)$.
\medskip

\subsection{Bessel functions, kernels}
The entire Bessel functions are
\begin{equation}\label{B}j_\alpha(z)=\Gamma(\alpha+1)\left(\frac{2}{z}\right)^\alpha J_\alpha(z)=\sum_{k=0}^\infty\frac{(-1)^k\Gamma(\alpha+1)}{\Gamma(k+1)\Gamma(k+\alpha+1)}\left(\frac{z}{2}\right)^{2k}.\end{equation}
Subsequently we assume that $\alpha>-\frac{1}{2}$. We note that $j_{-\frac{1}{2}}(x)=\cos x$.
\begin{equation}\label{j1}\|j_\alpha\|_{\infty, \mathbb{R}_+}=j_\alpha(0)=1.\end{equation}
\begin{equation}\label{j2}j_\alpha'(z)=-\frac{1}{2(\alpha+1)}zj_{\alpha+1}(z),\end{equation}
see \cite{be}.

We introduce the following abbreviation, cf. \cite{ss}.
$$\mathbb{j}_a(x,\xi):=\prod_{i=1}^nj_{\alpha_i}(x_i\xi_i).$$

The Bessel kernel is
\begin{equation}\label{g1}G_{a,\nu}(x):=\frac{2^{\frac{n-a-\nu}{2}+1}}{\Gamma\left(\frac{\nu}{2}\right)\prod_{i=1}^n\Gamma(\alpha_i+1)}\frac{K_{\frac{n+|a|-\nu}{2}}(|x|)}{|x|^{\frac{n+|a|-\nu}{2}}},\end{equation}
where $K_\alpha$ is the modified Bessel function of the second kind, see \cite{be}.

\subsection{Translations}

If $|a|>0$, on the positive orthant we introduce the generalized or Bessel translation defined by Delsarte (see \cite{d}).
The Bessel translation of a function, $f$ (see e.g. \cite{le}, \cite{p}, \cite{ss}) is
$$T_{a}^tf(x)=T_{a_n}^{t_n}\dots T_{a_1}^{t_1}f(x_1,\dots ,x_n),$$
where
\begin{equation}\label{tra1}T_{a_i}^{t_i}f(x_1,\dots ,x_n)$$ $$=\frac{\Gamma(\alpha_i+1)}{\sqrt{\pi}\Gamma\left(\alpha_i+\frac{1}{2}\right)}\int_0^\pi f(x_1,\dots,\sqrt{x_i^2+t_i^2 -2x_it_i\cos\vartheta_i}, x_{i+1}, \dots, x_n)\sin\vartheta^{2\alpha_i}d\vartheta_i.\end{equation}
For $|a|>0$, then we have
$$T_{a}^tf(x)=T_{a}^xf(t),$$
\begin{equation}\label{tt}T_{a}^yf(x)T_{a}^zf(x)=T_{a}^zf(x)T_{a}^yf(x),\end{equation}
see \cite[(7.1)]{le}.
If $|a|=0$, on $\mathbb{R}^n$ we use the standard translation, i.e.
$$T^t_{0}f(x):=f(x-t).$$
We have to note, that the Bessel translation with respect to $\alpha = -\frac{1}{2}$, i.e. $2\alpha+1=0$, is $T^t_{\dot{0}}f(x)=\frac{1}{2}(f(x+t)+f(|x-t|))$.\\
For all $a$ ($|a|\ge 0$)$, T_a$ is positive operator, and
\begin{equation}\label{t1}\|T_{a}^tf(x)\|_{p,a}\le \|f\|_{p,a}, \ws \ws 1\le p \le\infty,\end{equation}
for the Bessel case, see e.g. \cite{le}.

Of course, translations induce convolutions
$$f*_a g(t) =\int_{\mathbb{R}^n_{(+)}}T_{a}^tf(x)g(x)x^adx,$$
The standard convolution is denoted by $*=*_0$. In both cases we have
$$f*_a g =g*_a f,$$
\begin{equation}\label{ass}(f*_a g)*_a h =f*_a (g*_a h),\end{equation}
$f,g,h \in\mathcal{S}_{(e)}$, see e.g. \cite[(2.30)]{p}, and for $\mathcal{S}_{e}'$, see \cite{k}.\\
In the mixed space we introduce the next (mixed) convolution, c.f. [(3.174)]\cite{ss}.
\begin{equation}f*_Ag(x,t):=\int_0^t\int_{\mathbb{R}^n_{(+)}}T^y_af(x,\tau)g(y,t-\tau)y^adyd\tau\end{equation} $$=\int_0^\infty\int_{\mathbb{R}^n_{(+)}}T^y_af(x,\tau)g_+(y,t-\tau)y^adyd\tau=\int_0^\infty\int_{\mathbb{R}^n_{(+)}}T^y_af(x,\tau)T^\tau_{0}g_+(y,t)y^adyd\tau$$ $$=\int_{\mathbb{R}^{n+1}_{(+)}}T^\tau_{0}T^y_af(x,t)g(y,\tau)y^adyd\tau.$$
\begin{equation}\label{plus} g_+(u,v)=0 \ws \ws  \mbox{if} \ws v<0.\end{equation}
For $u,v\in \mathcal{S}_{(e),t}'$, as usual, $\langle u*_A v, \varphi\rangle_A= v\left((y,s)\mapsto u\left((x,\tau)\mapsto T^y\varphi(x, s+\tau)\right)\right),$
in particular $\delta_A*_Au=u$.

\subsection{Transformations}
In the sequel we use Fourier transformation on $\mathbb{R}^n$ if $|a|=0$, and Hankel transformation on $\mathbb{R}^n_+$ if $|a|>0$ with normalization as follows.
$$ \mathcal{I}_0f(y):=\mathcal{F}f(y)=\hat{f}(y)=\int_{\mathbb{R}^n}f(x)e^{-i\langle x,y\rangle}dx;$$ $$ \mathcal{I}_af(y):=\mathcal{H}_a(f,\xi)=\tilde{f}(\xi)=\int_{\mathbb{R}^n_+}f(x)\mathbb{j}_a(x,\xi)x^adx.$$
If $f\in L^1_a(\mathbb{R}^n_+)$ $(|a|>0)$ is of bounded variation in a neighborhood of a point $x$ of continuity of $f$, then the next inversion formula holds, see [page 38.]\cite{ss}.
$$f(x)=\mathcal{H}_a^{-1}(\tilde{f}(\xi))(x)=\frac{2^{n-|a|}}{\prod_{i=1}^n\Gamma^2(\alpha_i+1)}\int_{\mathbb{R}^n_+}\tilde{f}(\xi)\mathbb{j}_a(x,\xi)\xi^ad\xi.$$
In both cases $\mathcal{I}_a$ maps $L^p_a$ to $L^{p'}_a$, $1\le p\le 2$. (For Hankel case see e.g. \cite{ss}.)\\
The relation to the appropriate convolutions is
\begin{equation}\label{hkon}\mathcal{I}_a(f*_ag)=\mathcal{I}_a(f)\mathcal{I}_a(g),\end{equation}
see e.g \cite[(3.168)]{ss} (Hankel case).

\medskip

\note
$S^n_1\subset \mathbb{R}^n$ is the ($n-1$-dimensional) unit sphere, i.e. $\{\Theta\in \mathbb{R}^n:|\Theta|=1\}$. $S_{1+}^n=S_{1}^n\cap \mathbb{R}^n_+$.
\begin{equation}\label{felsz0}|S_{1}^n|_0=\int_{S_{1}^n}1dS(\Theta)=\frac{2\pi^{\frac{n}{2}}}{\Gamma\left(\frac{n}{2}\right)},\end{equation} and
\begin{equation}\label{felsz}|S_{1+}^n|_a=\int_{S_{1+}}\Theta^adS(\Theta)=\frac{\prod_{i=1}^n \Gamma(\alpha_i+1)}{2^{n-1}\Gamma\left(\frac{n+|a|}{2}\right)},\end{equation}
see \cite[(2.8)]{esk}.

Subsequently, the following formulae will be of key importance, stating that both Fourier and the Hankel transforms of a radial function can be expressed as a one-dimensional Hankel transform.

\medskip

\begin{lemma}\label{le1} Let $f\in L^1_a(\mathbb{R}^n_{(+)})\cap L_a^2(\mathbb{R}_{(+)}^n)$ be a radial function, i.e. $f(x)=\varphi(|x|)$, where $\varphi$ is a function of one variable. Then denoting by $r=|x|$\\
for $|a|=0$, i.e. Fourier transform case, we have
\begin{equation}\label{rf}\mathcal{F}f(|\xi|)= |S_{1}^n|\int_0^\infty \varphi(r)j_{\frac{n}{2}-1}(r |\xi|)r^{n-1}dr= |S_{1}^n|\mathcal{H}_{\frac{n}{2}-1}(\varphi)(|\xi|),\end{equation}
see \cite[Theorems 2,3]{o}\label{lor}.\\
In $|a|>0$ case the Hankel transform of a radial function,  $\mathcal{H}_a(f)$, is also radial, furthermore we have
\begin{equation}\label{rh}\mathcal{H}_a(f)(\xi)=|S_{1+}^n|_a\int_0^\infty \varphi(r)j_{\frac{n+|a|}{2}-1}(|\xi|r)r^{n+|a|-1}dr,\end{equation}
see \cite[Lemma 3.2]{esk}.
\end{lemma}

\medskip

The Fourier-Hankel transform of $f\in L^1_A(\mathbb{R}^{n+1}_+)$
$$\mathcal{I}_A(f,\xi,\tau)=\int_0^\infty\int_{\mathbb{R}^n_{(+)}}f(x,t)\mathbb{j}_a(x,\xi)x^adxe^{-it\tau} dt.$$
(It should be mentioned here that there is different transformation in the literature called the "Fourier-Bessel transformation".)\\
Recalling that the Bessel functions are even, for $f\in \mathcal{S}_{(e)}'$ we have
\begin{equation}\label{HT} \mathcal{I}_a(T^y_af)(\xi)= \Phi_a(\xi,y)\mathcal{I}_a(f)(\xi); \ws \ws \ws T^y_a(\mathcal{I}_a f)(\xi)=\mathcal{I}_a(\Phi_a(x,-y)f(x)),\end{equation}
where
\begin{equation}\Phi_a(u,y)=\mathbb{j}_a(u,y), \ws |a|>0; \ws \ws \Phi_a(u,y)=e^{-i\langle u,y\rangle}, \ws a=0,\end{equation}
see e.g. \cite[(3.158)]{ss} for the Bessel case.\\

For $u\in \mathcal{S}_{(e)}'$ and $\varphi\in \mathcal{S}_{(e)}$ we have (see e.g. \cite[(1.83)]{ss})
$$\langle \mathcal{I}_au, \varphi \rangle_a=\langle u, \mathcal{I}_a\varphi \rangle_a.$$

In particular if $\varphi \in \mathcal{S}_{(e)}$,
$$\langle \mathcal{H}_a\delta_a, \varphi \rangle_a=\langle \delta_a, \mathcal{H}_a\varphi \rangle_a=\tilde{\varphi}(0)=\int_{\mathbb{R}^n_+}\varphi(x)\mathbb{j}_a(x,0)x^adx=\langle 1, \varphi \rangle_a,$$
and similarly if $\varphi \in \mathcal{S}_{(e),t}$,
$$\langle \mathcal{I}_A\delta_A, \varphi \rangle_A=\mathcal{I}_A(\varphi)(0,0)=\langle 1, \varphi \rangle_A,$$
That is
\begin{equation}\label{deltatraf}\mathcal{I}_a(\delta_a) = \mathcal{I}_A\delta_A =1.\end{equation}

\subsection{Fractional heat equations}

If $\alpha>-\frac{1}{2}$,
$$B_{\alpha}:= B_{\alpha,x} =\frac{\partial^2}{\partial x^2} +\frac{2\alpha+1}{x}\frac{\partial}{\partial x}.$$
\begin{equation}\label{sfv}B_{\alpha}j_\alpha(\lambda x)=-\lambda^2j_\alpha(\lambda x),\end{equation}
see [(2.11)]\cite{le}
The Bessel-Laplace operator (in $n$-dimension) is defined as
$$\Delta_a:=\Delta_{a,x}=\sum_{i=1}^nB_{\alpha_i,x_i},$$
where $x\in \mathbb{R}^n_+$.\\
If $a=0$ we arrive to the standard Laplace operator,
$$\Delta_0=\Delta.$$
For $f\in\mathcal{S}_{(e)}'$ we have
\begin{equation}\label{TB}T^y(\Delta_a f)=\Delta_a (T^yf),\end{equation}
see e.g. \cite[(2.20)]{p}, and
\begin{equation}\label{HB}  \mathcal{I}_a(-\Delta_a f)=|\xi|^2  (\mathcal{I}_af)(\xi),\end{equation}
see \cite[(1.95)]{ss}.

The fractional Laplace and Laplace-Bessel operators are defined via Fourier or Hankel transforms as it follows, see \cite{bg},\cite{l}. Let $0< \gamma <2$ and $f, |\xi|^{\gamma}(\mathcal{I}_af)(\xi)\in L^2_a $. Then
\begin{equation}\label{frach} \mathcal{I}_a\left((-\Delta_a)^{\frac{\gamma}{2}}f\right)(\xi)=|\xi|^{\gamma}(\mathcal{I}_a f)(\xi).\end{equation}
The other way to define the fractional operators is
\begin{equation}\label{fracl} (-\Delta_a)^{\frac{\gamma}{2}}f(x)=C(n,a,\gamma)\int_{\mathbb{R}^n_{(+)}} \frac{f(x)-T^y_af(x)}{|y|^{n+|a|+\gamma}}y^ady, \end{equation}
where in the second case the integral is meant in principal value sense. The Bessel case see \cite{l} and in one dimension \cite{bg}.

Define the next fractional operators
\begin{equation}\label{bfo} L_{\gamma,a}:=(-\Delta_a)^{\frac{\gamma}{2}}+\partial_t. \end{equation}
By standard arguments the fundamental solution is given by the corresponding transform:
\begin{equation}\label{E} P_{\gamma,a}(x,t)=\left\{\begin{array}{ll}\mathcal{I}_a^{-1}\left(e^{-{|\xi|^\gamma}t} \right), \ws \mbox{if}  \ws t>0,\\0, \ws \mbox{if}  \ws t\le 0,\end{array}\right.\end{equation}
where the inverse Fourier and Hankel transforms refer to $\xi$ and $x\in \mathbb{R}^n_{(+)}$, c.f. [Theorem 4.5]\cite{bg}, see also \cite{v}.

\medskip

In view of Lemma \ref{le1} we can investigate Laplace and Bessel cases at the same time.
That is $P_{\gamma,a}(x,t)=P_{\gamma,a}(|x|,t)$ is radial function of $x$, moreover denoting by $|\xi|=\varrho$ we have
\begin{equation}\label{Pm} P_{\gamma,a}(x,t)=c(n,a)\int_0^\infty e^{-t\varrho^\gamma}j_{\frac{n+|a|-2}{2}}(\varrho |x|)\varrho^{n+|a|-1}d\varrho,\end{equation}
where

$$c(n,a)=\left\{\begin{array}{ll}\frac{1}{2^{|a|-1}\prod_{i=1}^n\Gamma(\alpha_i+1)\Gamma\left(\frac{n+|a|}{2}\right)}, \ws |a|>0\\
\frac{|S_1^n|_0}{(2\pi)^n}, \ws a=0. \end{array}\right.$$

\medskip

In view of \cite[ (4.1)-(4.3)]{bg}, $x\in \mathbb{R}_+$
\begin{equation}\label{33} E_{\gamma,\nu}(x,t):=\frac{1}{2^{\nu+1}\Gamma(\nu+1)}\int_0^\infty  e^{-t\varrho^\gamma}j_\nu(\varrho x)\varrho^{2\nu+1}d\varrho =t^{-\frac{2(\nu+1)}{\gamma}}k_{\gamma,\nu}\left(t^{-\frac{1}{\gamma}}x\right).\end{equation}
Let $\nu=\frac{n+|a|-2}{2}$. Then $\nu\ge -\frac{1}{2}$ and for $t>0$ we have
\begin{equation}\label{pp}P_{\gamma,a}(x,t)=c(n,a)E_{\gamma,\frac{n+|a|}{2}-1}(|x|,t),\end{equation}
where
\begin{equation}c(n,a)=\left\{\begin{array}{ll}\frac{2^{\frac{n-|a|}{2}+1}}{\prod_{i=1}^n\Gamma(\alpha_i+1)}, \ws |a|>0\\
\frac{2}{(2\pi)^{\frac{n}{2}}}, \ws a=0. \end{array}\right.\end{equation}

The expansion of $k_{\gamma,\nu}$ is given in \cite[Proposition 4.1]{bg}. Together with \eqref{33} and \eqref{pp} it gives

\medskip

\begin{lemma}\label{lk} $x,t\neq 0$. Then we have the following expansions.
If $0<\gamma<1$,
$$P_{\gamma,a}(x,t)=c(n,a)\frac{2^{n+1}}{\prod_{i=1}^n\Gamma(\alpha_i+1)}\frac{t}{|x|^{n+|a|+\gamma}}$$ \begin{equation}\label{mind}\times\sum_{l=0}^\infty (-1)^{l}\frac{\Gamma\left(\frac{(l+1)\gamma}{2}+1\right)\Gamma\left(\frac{n+|a|+(l+1)\gamma}{2}\right)}{(l+1)!}\sin\left(\frac{\gamma (l+1) \pi}{2}\right)2^{\gamma (l+1)}\frac{t^l}{|x|^{\gamma l}}.\end{equation}
\begin{equation}P_{1,a}(x,t)=c(n,a)\frac{t}{(t^2+|x|^2)^{\frac{n+|a|+1}{2}}},
\end{equation}
If $1<\gamma<2$,
\begin{equation}\label{nagy}P_{\gamma,a}(x,t)=c(n,a)\frac{1}{2^{|a|}\gamma \prod_{i=1}^n\Gamma(\alpha_i+1)}\frac{1}{t^{\frac{n+|a|}{\gamma}}}\sum_{l=0}^\infty \frac{(-1)^{l}}{l!}\frac{\Gamma\left(\frac{2l+n+|a|)}{\gamma}\right)}{\Gamma\left(l+\frac{n+|a|}{2}\right)4^l}\frac{|x|^{2l}}{t^{\frac{2l}{\gamma}}}.\end{equation}
\end{lemma}

\medskip

\begin{lemma}\label{E}\cite[Proposition 4.2]{bg}
$$\frac{1}{2^{\nu+1}\Gamma(\nu+1)}\int_0^\infty E_{\gamma,\nu}(x,t)x^{2\nu+1}dx=1.$$
\end{lemma}

\medskip

\noindent{\bf Remark.}\\

\noindent {\bf 1.} $P_{\gamma,a}(x,t)>0$, see \cite{bg}.\\

\noindent {\bf 2.}\\
In view of \eqref{j2} and \eqref{Pm} we have
\begin{equation}\label{dxP}\partial_{x_i}P_{\gamma,a}(x,t)=-c(n,a)x_iP_{\gamma,^ia}(x,t),\end{equation}
where $c(n,a)$ is a positive constant, and $^ia:=a_1,\dots , a_{i-1}, a_i+2, a_{i+1}, \dots , a_n$. \\

\noindent {\bf 3.}\\
Le $b\in \mathbb{R}^n_+$ be a multiindex. Then, by Lemma \ref{E} we have
\begin{equation}\label{p1}\int_{\mathbb{R}^n_+}P_{\gamma,a+b}(x,t)|x|^{|b|}x^adx=c(n,a,b)\int_{\mathbb{R}^n_+}P_{\gamma,a+b}(x,t)x^{a+b}dx<\infty.\end{equation}
Indeed,
$$\int_{\mathbb{R}^n_+}P_{\gamma,a+b}(x,t)|x|^{|b|}x^adx$$ $$=c(n,a,b)|S_{1+}|_a\int_0^\infty E_{\gamma,\frac{n+|a|+|b|}{2}-1}(r,t)r^{n+|a|+|b|-1}dr.$$

\medskip

\subsection{Capacity and Sobolev-type spaces}

In this subsection we introduce some Sobolev-type spaces and different capacities.

\medskip

\note

Let $1\le p<\infty$. Either $|a|>0$ or $a=0$, we define the Sobolev-type spaces $W^{m,p}_{\Delta_a}$ and $L^p_{a,\nu}$ as follows (c.f. \cite{h}, \cite[(1.9)]{p} and \cite[(1.2.29)]{ah}).
$$W^{m,p}_{\Delta_a}:=\left\{f: \mathbb{R}^n_{(+)} \to \mathbb{R} : \Delta_a^k f \in  L^p_a, \ws k=0, \dots , m \right\},$$ $$ \|f\|_{W^{m,p}_{\Delta_a}}=\left(\sum_{k=0}^m \|\Delta_a^k f\|_{p,a}^p\right)^{\frac{1}{p}}.$$
$$L^p_{a,\nu}:=L^p_{a,\nu}(\mathbb{R}^n_+)=\{f :f=G_{a,\nu}*_a g; \ws   g\in L^p_a\}, \ws \ws \ws \|f\|_{p,a,\nu}:=\|g\|_{p,a}.$$
Moreover $f(x,t)\in W^{m,p}_{\Delta_a,A}$ if $\|f_t\|_{W^{m,p}_{\Delta_a}}\in L^1(\mathbb{R}_+)$, where $f_t(x)$ stands for $f(x,t)$.

It is clear that the above defined $\Delta_a$-Sobolev spaces are Banach spaces; $W^{m,p}\subsetneq W^{m,p}_{\Delta}$, and like $W^{m,p}$ these are separable if $1\le p<\infty$ and are reflexive and uniformly convex if $1<p<\infty$, cf. \cite[3.4]{a}.\\
Since  these do not coincide with  $W^{m,p}$, Calder\'on's Theorem (see e.g. \cite[Theorem 1.2.3]{ah}) is not applicable. Instead of this, in Bessel case we have the following lemma.

\medskip

\begin{lemma}\label{l1}\cite[Lemma 2.]{h} Let $|a|>0$ and $m$ be a positive integer, $1<p<\infty$. Then (with equivalent norms)
$$W^{m,p}_{\Delta_a}=L^p_{a,2m}.$$
\end{lemma}

\medskip

Let $1\le p<\infty$. We define the next dual $\Delta_a$-Sobolev spaces.

\medskip

\note

As in the standard case let $W^{m,p}_{0,\Delta_a}$ be the closure of $C_0^\infty(\mathbb{R}^n_{(+)})$ in $W^{m,p}_{\Delta_a}$, and $W^{-m,p'}_{\Delta_a}$ stands for its dual space. Moreover  $W^{-m,p'}_{\Delta_a,A}$ is the dual of $W^{m,p}_{0,\Delta_a,A}$.  Notice that $\|\cdot\|_{W^{m,p}_{\Delta}}\le \|\cdot\|_{W^{m,p}}$, thus $W^{m,p}_{0}\subset W^{m,p}_{0,\Delta}$.  The dual space possesses the following characterization.

\medskip

\begin{proposition}\label{pdu} $W^{-m,p'}_{\Delta_a}$ is isometrically isomorph with the space of those distributions $T\in \mathcal{D}'(\mathbb{R}^n_{(+)})$ which have the form $T=\sum_{k=0}^m \Delta_a^kT_{v_k}, \ws v_k\in L^{p'}_a$ and $\|T\|=\inf\left\{\left(\sum_{k=0}^m\|v_k\|_{p',a}^{p'}\right)^{\frac{1}{p'}}: T=\sum_{k=0}^m \Delta_a^kT_{v_k}; \ws \right\}$, where $T_{v_k}$ is the regular distribution generated by $v_k\in L^{p'}_a$.
\end{proposition}

The proof of Proposition \ref{pdu} follows the standard arguments,see \cite[Chapter III]{a}, so we omit the details.

\medskip

\note

Recalling \eqref{g1}, let
$$\nu_G:=\nu_t*_aG_{a,2},$$

In the spirit of "dual definition of capacity", see \cite[Theorem 2.7.2]{ah}, with this notation we introduce different capacities.

\begin{defi}\label{dcap} Let $K\subset \mathbb{R}^{n+1}_{(+)}$ be compact.
\begin{equation}\label{dncap} N_{a,\gamma,p'}^{\frac{1}{p'}}(K):=\sup\{\langle \nu,1\rangle_A : \nu\in \mathcal{S}_{(e),t}', \ws \mathrm{supp}\nu \subset K, \ws \|P_{\gamma,a}*_A\nu_G\|_{p,A}\le 1\}.\end{equation}

\begin{equation}\label{ncap} Z_{a,\gamma,p'}^{\frac{1}{p'}}(K):=\sup\{\langle \nu,1\rangle_A : \nu\in \mathcal{S}_{(e),t}', \ws \mathrm{supp}\nu \subset K, \ws \|P_{\gamma,a}*_A\nu\|_{W^{-1,p}_{\Delta_a},A}\le 1\}.\end{equation}

\begin{equation}\label{dcap} C_{a,\gamma,p'}^{\frac{1}{p'}}(K):=\sup\{\langle \nu,1\rangle_A : \nu\in \mathcal{S}_{(e),t}', \ws \mathrm{supp}\nu \subset K, \ws \|P_{\gamma,a}*_A\nu\|_{p,A}\le 1\}.\end{equation}
\end{defi}

Since $G_{a,2}\in L^1_a$ it is clear that $C_{a,\gamma,p'}\le N_{a,\gamma,p'}$, and also  $C_{a,\gamma,p'}\le Z_{a,\gamma,p'}$.

\section{Removable sets}

We say that a compact set $K\subset \mathbb{R}^{n+1}_{(+)}$ is removable for $L_{\gamma,a}$ in $L^p_A$ if for any $f\in L^p_A$, $L_{\gamma,a}f=0$ on $\mathbb{R}^{n+1}_{(+)}\setminus K$, then $f$ also satisfies the fractional heat diffusion equation on $\mathbb{R}^{n+1}_{(+)}$.

\begin{proposition}\label{lp} Let $1< \gamma<2$ and $1< p <\infty $. Then either for $|a|>0$ or $a=0$, if $Z_{a,\gamma,p'}^{\frac{1}{p'}}(K)=0$, then $K\subset \mathbb{R}^{n+1}_{(+)}$ is removable for $L_{\gamma,a}$ in $L^p_A$.
\end{proposition}

In \cite{p} the equivalence of Bessel $p$-modulus of smoothness and the corresponding K-functional is proved in one-dimensional case. To prove Proposition \ref{lp} we need something similar. To this we define the next difference.

\begin{defi}\label{sd}
\begin{equation}\label{bedi} \Delta_{r,s,a}f(x):=\frac{1}{|S^n_{1(+)}|_a}\int_{S^n_{1(+)}}(f(x)-T^{r\Theta}_af(x))\Theta^adS(\Theta),\end{equation}
where
$$S^n_{1(+)}=\left\{\begin{array}{ll}S_{1}^n, \ws a=0\\ S_{1+}^n, \ws |a|>0;\end{array}\right. $$
cf. \eqref{felsz} and $dS$ stands for the integration on the sphere.\end{defi}

Subsequently we need the next formulae.

\medskip

\begin{lemma}\begin{equation}\label{besatl} \frac{1}{|S_{1+}^n|_a}\int_{S_{1+}^n}\mathbb{j}_a(r\Theta,\xi))\Theta^adS(\Theta)=j_{\frac{n+|a|}{2}-1}(r|\xi|),\end{equation}
c.f. \cite[(3.140)]{ss}.

\begin{equation}\label{fours}\frac{1}{|S_{1}^n|}\int_{S_{1}^n}e^{-ir|\xi|\langle \Theta,\vartheta\rangle}dS(\Theta)=j_{\frac{n}{2}-1}(r|\xi|),
\end{equation}
c.f. \cite[page 136]{sch} and the references therein.
\end{lemma}

\medskip

Furthermore we have the following lemma.

\medskip

\begin{lemma}\cite[Lemma 4.1]{p}\label{lp41}
$$\mathcal{H}_\alpha^{-1}\left(\frac{1-j_{\alpha}(u)}{u^2}\right)\in L^1_\alpha(\mathbb{R}_+).$$
\end{lemma}

\medskip

\begin{lemma}\label{gamma1} Let $1\le p \le \infty$.
 If $f\in \mathcal{S}_{(e)}$, then
\begin{equation}\label{mK} \|\Delta_{r,s,a}f(x)\|_{p,a}\le c r^2\|\Delta_a f\|_{p,a},\end{equation}
respectively, and $c$ is independent of $f$ and $r$.
\end{lemma}

\proof
In view of \eqref{HT} and \eqref{besatl}, the translation rule of Fourier transform and by \eqref{fours} we have
$$\mathcal {I}_a(\Delta_{r,s,a}f)(\xi)=r^2\frac{1-j_{\frac{n+|a|}{2}-1}(r|\xi|)}{|\xi|^2r^2}|\xi|^2(\mathcal{I}_a f)(\xi).$$
Then by the inversion formulae, \eqref{HB} and \eqref{hkon} we have
$$\|\Delta_{r,s,a}f\|_{p,a}\le r^2\left\|\mathcal{I}_a^{-1}\left(\frac{1-j_{\frac{n+|a|}{2}-1}(r|\xi|)}{|\xi|^2r^2}\right)\right\|_{1,a}\|(-\Delta_a)f\|_{p,a},$$
Lemma \ref{le1} and a simple replacement imply
$$\left\|\mathcal{I}_a^{-1}\left(\frac{1-j_{\frac{n+|a|}{2}-1}(r|\xi|)}{|\xi|^2r^2}\right)\right\|_{1,a}$$ $$=c(n,a)\left\|\mathcal{H}_{\frac{n+|a|}{2}-1}^{-1}\left(\frac{1-j_{\frac{n+|a|}{2}-1}(\varrho)}{\varrho^2}\right)\right\|_{1,n+|a|-1},$$
where $\varrho=|\xi|$. Finally Lemma \ref{lp41} ensures the result.

\medskip

With the introduced spherical difference we estimate for the norm of the intermediate fractional derivatives by the $\Delta_a$-Sobolev norm of the function. This is similar to the estimate for the so-called Gagliardo seminorm, cf. \cite[Proposition 2.2]{dpv}.

\medskip

\begin{lemma}\label{gammadelta}
Let $\gamma \in (0,2)$, $1\le p<\infty$. If $f\in \mathcal{S}_{(e)}$, then
\begin{equation}\label{D}\left\|(-\Delta_a)^{\frac{\gamma}{2}}f\right\|_{p,a}\le c \|f\|_{W^{1,p}_{\Delta_a}},\end{equation}
where $c=c(n,a,p)$.
\end{lemma}

\proof For any $\eta>0$
$$\left\|(-\Delta_a)^{\frac{\gamma}{2}}f\right\|_{p,a}$$ $$\le \left\|\int_{|y|<\eta}\frac{T^y_af(x)-f(x)}{|y|^{n+|a|+\gamma}}y^ady\right\|_{p,a}+\left\|\int_{|y|\ge\eta}\frac{T^y_af(x)-f(x)}{|y|^{n+|a|+\gamma}}y^ady\right\|_{p,a}=I+II.$$
By Minkowski inequality and then considering Lemma \ref{gamma1} we have
$$I=\left(\int_{\mathbb{R}^n_{(+)}}\left|\int_0^\eta \frac{1}{r^{1+\gamma}}\int_{S_{1(+)}}\left(T^{r\Theta}_af(x)-f(x)\right)\Theta^adSdr\right|^px^adx\right)^{\frac{1}{p}}$$ $$\le |S_{1(+)}|_a\int_0^\eta \frac{1}{r^{1+\gamma}}\|\Delta_{r,s,a}f(x)\|_{p,a}dr$$
$$\le c \|\Delta_a f\|_{p,a}\int_0^\eta \frac{1}{r^{\gamma-1}}dr=c(n,a,\gamma,\eta)\|\Delta_a f\|_{p,a}.$$
Again by Minkowski inequality and by \eqref{t1} we have
$$II\le 2\|f\|_{p,a}\int_{|y|\ge\eta}\frac{1}{|y|^{n+|a|+\gamma}}y^ady= c(n,a,\gamma,\eta)\|f\|_{p,a}.$$
Choosing $\eta=1$, say the lemma is proved.

\medskip

\proof (of Proposition \ref{lp})\\
Let $K\subset \mathbb{R}^n_{(+)}$ be a compact set and $Q\subset \mathbb{R}^n_{(+)}$ be a cube so that $K\subset Q$. $Q= Q_s\times I_T$. Indirectly, assume that $K$ is not removable. Then there exists an $f\in L^p_A$ such that $L_{\gamma,a}f\equiv 0$ on $Q\setminus K$, and $L_{\gamma,a}f\not\equiv 0$ on $Q$. We may assume that $f$ is supported on $Q$. Define the distribution $\nu:=L_{\gamma,a}f$. Since the support of $\nu$ is in $Q$ and $\nu \not\equiv 0$ on $Q$, there is a $\varphi \in C_0^\infty$, supported in $Q$, such that
$$\langle\nu,\varphi\rangle_A=\langle\nu,\varphi\rangle_A=\langle \varphi\nu,1\rangle_A>0.$$
In view of Definition \ref{dcap} it is enough to show that $\|P_{\gamma,a}*_A\varphi\nu\|_{W^{-1,p'}_{\Delta_a},A}<\infty .$\\
Let $\Psi\in \mathcal{S}_{(e),t}$. General Parseval's formula ensures that
$$\langle P_{\gamma,a}*_A\varphi\nu,\Psi\rangle_A=\langle \nu, \varphi\left(P_{\gamma,a}*_A\Psi\right)\rangle_A= \langle \nu, P_{\gamma,a}*_A L_{\gamma,a}\left(\varphi\left(P_{\gamma,a}*_A\Psi\right)\right)\rangle_A$$ $$=\langle P_{\gamma,a}*_A\nu, L_{\gamma,a}\left(\varphi\left(P_{\gamma,a}*_A\Psi\right)\right)\rangle_A.$$
Since $P_{\gamma,a}$ is the fundamental solution, $P_{\gamma,a}*_A\nu \in L^p_a$, that is
$$|\langle P_{\gamma,a}*_A\varphi\nu,\Psi\rangle_A|\le \|P_{\gamma,a}*_A\nu\|_{p,A}\|L_{\gamma,a}\left(\varphi\left(P_{\gamma,a}*_A\Psi\right)\right)\|_{p',A}.$$
Since $K\subset Q_s\subset \mathbb{R}^n_{(+)}$, it is readily seen that for all $t\in I_T$ $\left(\varphi\left(P_{\gamma,a}*_A\Psi\right)\right)_t\in \mathcal{S}_{(e)}$. So by Lemma \ref{gammadelta}
\begin{equation}\label{eee}\|L_{\gamma,a}\left(\varphi\left(P_{\gamma,a}*_A\Psi\right)\right)\|_{p',A}\le c\left(\left(\int_{I_T}\left\|\Delta_a\left(\varphi\left(P_{\gamma,a}*_A\Psi\right)\right)(x,t)\right\|_{p',a}^{p'}dt\right)^{\frac{1}{p'}}\right.\end{equation} $$\left.+\left(\int_{I_T}\left\|(\varphi\left(P_{\gamma,a}*_A\Psi\right)(x,t)\right\|_{p',a}^{p'}dt\right)^{\frac{1}{p'}}+\left(\int_{I_T}\left\|\partial_t\left(\varphi\left(P_{\gamma,a}*_A\Psi\right)\right)(x,t)\right\|_{p',a}^{p'}dt\right)^{\frac{1}{p'}}\right)$$
$$\le c\left(\|(\Delta_a \varphi)(P_{\gamma,a}*_A\Psi)\|_{p',A,Q}+\left\|\sum_{i=1}^n\partial_{x_i} \varphi \partial_{x_i}(P_{\gamma,a}*_A\Psi)\right\|_{p',A,Q}\right.$$ $$\left.+\|\varphi\Delta_a (P_{\gamma,a}*_A\Psi)\|_{p',A,Q}+\| \varphi(P_{\gamma,a}*_A\Psi)\|_{p',A,Q}+\|(\partial_t \varphi)(P_{\gamma,a}*_A\Psi)\|_{p',A,Q}\right.$$ $$\left.+\|\varphi \partial_t(P_{\gamma,a}*_A\Psi)\|_{p',A,Q}\right)=c(I+II+III+IV+V+VI).$$
Let $t_0$ be such that $I_T\subset (0, t_0)=:I_{t_0}$.  Let $h\in \mathcal{S}_{(e),t}$. Minkowski inequality, \eqref{t1} and \eqref{p1} imply
$$\|P_{\gamma,a}*_Ah\|_{p',A,Q}\le \int_{I_{t_0}}\int_{\mathbb{R}^n_{(+)}}P_{\gamma,a}(y,\tau)\|T^yh(x,t-\tau)\|_{p',A}y^adyd\tau\le c t_0 \|h\|_{p',A}.$$
Thus, considering that \eqref{TB} implies that $\Delta_a (P_{\gamma,a}*_A\Psi)= P_{\gamma,a}*_A\Delta_a\Psi)$, we have
$$I\leq c t_0\|\Delta_a \varphi\|_\infty  \|\Psi\|_{p',A}; \ws \ws III\le  c t_0\|\varphi\|_\infty \|\Delta_a\Psi\|_{p',A};$$ $$ IV\le c t_0\|\varphi\|_\infty \|\Psi\|_{p',A}; \ws \ws V\le c t_0\|\partial_t\varphi\|_\infty \|\Psi\|_{p',A}.$$

To estimate the remainder terms we compute the derivatives of the kernel.
Note here that the order of the first derivative and the Bessel translation cannot be interchanged. In view of \eqref{dxP} in Bessel case we have
$$c(n,a)\left\|\partial_{x_i} \varphi \partial_{x_i}(P_{\gamma,a}*_A\Psi)\right\|_{p',A,Q}$$ $$\le \|\partial_{x_i} \varphi\|_\infty \left(\int_{Q}\left|\int_0^t\int_{\mathbb{R}^n_{(+)}}\Psi(y,t-\tau)\int_{([0,\pi]^n}(x_i-y_i\cos\vartheta_i)\right.\right.$$ $$ \left.\left. \times E_{\gamma,\frac{n+|a|}{2}}\left(\sqrt{\sum_{j=1}^n(x_j^2+y_j^2-2x_jy_j\cos\vartheta_j)},\tau\right))d\sigma(a,\vartheta)y^ady\right|^{p'}x^adxdt\right)^{p'}$$ $$=\|\partial_{x_i} \varphi\|_\infty J_i.$$
By the standard replacement $x_j-y_j\cos\vartheta_j=z_j\cos\varphi_j$, $y_j\sin\vartheta_j=z_j\sin\varphi_j$, $j=1,\dots ,n$ the translation is transferred to $\Psi$ and by Minkowski inequality, as above, we have
$$J_i=\left\|\int_0^t\int_{\mathbb{R}^n_+}\int_{([0,\pi]^n}\Psi(\dots,\sqrt{x_j^2+z_j^2-2x_jz_j\cos\varphi_j},\dots ;t-\tau)\right.$$ $$\left. \times z_i\cos\varphi_i E_{\gamma,\frac{n+|a|}{2}}(|z|,\tau)d\sigma(a,\varphi)z^adz\right\|_{p',A,Q}$$ $$\le \int_0^{t_0}\int_{\mathbb{R}^n_+}z_i E_{\gamma,\frac{n+|a|}{2}}(|z|,\tau)\left\|T^x|\Psi|(z,t-\tau)\right\|_{p',A}z^adzd\tau$$ \begin{equation}\label{Ebecs}\le \|\Psi\|_{p',A}\left(\int_{B(0,1)_{n+1,+}}z_i E_{\gamma,\frac{n+|a|}{2}}(|z|,\tau)z^adzd\tau+\int_{(|(z,\tau)|>1)\cap (t<t_0)}(\cdot)z^adzd\tau\right)\end{equation} $$=\|\Psi\|_{p',A}(J_{i,1}+J_{i,2}).$$
First we estimate the integral of $z_i E_{\gamma,\frac{n+|a|}{2}}(|z|,\tau)$ on the positive orthant of the $n+1$ dimensional unit ball. According to Lemma \ref{lk}, on bounded sets
\begin{equation}\label{pkb}P_{\gamma,a}(x,t)\sim \frac{t}{\left(t^{\frac{2}{\gamma}}+|x|^2\right)^{\frac{n+|a|+\gamma}{2}}},\end{equation}
cf. also \cite[(2.4)]{v} if $a=0$.
Thus we introduce the modified spherical coordinates: $z_j=r\prod_{l=1}^{j-1}\sin\eta_l\cos\eta_j$, $j=1,\dots , n$; $\tau=\left(r\prod_{l=1}^{n}\sin\eta_l\right)^\gamma$ with Jacobian $J=\gamma r^{n-1+\gamma}\prod_{l=1}^n\sin^{n-1-l+\gamma}\eta_l$.
$$J_{i,1}=c(\gamma,n,a)\int_0^1r^{\gamma-2}dr,$$
For $\gamma>0$ the integral on the sphere is finite. Since $\tau$ is bounded, on exterior of the unit ball we can estimate the derivative of the kernel by a constant if $|z|$ is small. If $|z|>1$, say, there is a $j$ such that $z_i\le z_iz_j$, and in view of \eqref{p1}
$$J_{i,2}\le c t_0 \int_{\mathbb{R}^n_+}z_i z_jE_{\gamma,\frac{n+|a|}{2}}(|z|,\tau)z^adz<\infty.$$
Thus, if $1<\gamma<2$, $\partial_{x_i}P_{\gamma,a}(x,t)\in L^1_{A,\mathbb{R}^n_{(+)}\times I_{t_0}}$ and so in Bessel case $II$ is estimated as well.\\
In the standard case $\partial_{x_i}(P_{\gamma,0}*_A\Psi)=c\int_0^t\int_{\mathbb{R}^n}(x_i-y_i)E_{\gamma,\frac{n}{2}}(x-y,\tau)\Psi(x,t-\tau)dyd\tau$, and by the standard replacement we immediately arrive to \eqref{Ebecs}. Than we complete this part as in Bessel case.\\
To estimate $IV$ let us recall, that $P_{\gamma,a}$ is the fundamental solution. Thus
$$\|\partial_t(P_{\gamma,a}*_A\Psi)\|_{p',A,Q}\le \|\Psi\|_{p',A,Q}+\|(-\Delta_a)^{\frac{\gamma}{2}}(P_{\gamma,a}*_A\Psi)\|_{p',A,Q}=:D.$$
According to Lemma \ref{gammadelta}
$$D\le \|\Psi\|_{p',A}+c \left(\int_{I_t}\|\Delta_a(P_{\gamma,a}*_A)\Psi\|_{p',a}^{p'}dt\right)^{\frac{1}{p'}}\le c\|\Psi\|_{W^{1,p'}_{\Delta_a},A},$$
where we proceeded as in $III$.

\medskip

The curiosity of Proposition \ref{lp} is that the Laplace and Bessel-Laplace cases are treated simultaneously. Of course, this way the advantage of special properties are lost. The theorems below are formulated according to the specifics of the two different operators.

\medskip

\begin{theorem}\label{th1} Let $1< \gamma<2$ and $1< p <\infty $ and $|a|>0$. In this case if $N_{a,\gamma,p'}^{\frac{1}{p'}}(K)=0$, then $K\subset \mathbb{R}^{n+1}_+$ is removable for $L_{\gamma,a}$, and if $K\subset \mathbb{R}^{n+1}_+$ is removable for $L_{\gamma,a}$, then $C_{a,\gamma,p'}^{\frac{1}{p'}}(K)=0$.
\end{theorem}

\proof

According to Lemma \ref{l1} for any $t>0$ there is a $h_t(x):=h(x,t)\in \mathcal{S}_{e,t}$, such that $\Psi_t=G_{a,2}*_ah_t$ and $\|h_t\|_{p',a} \sim \|\Psi_t\|_{W^{1,p'}_{\Delta_a}}$.
$$\langle P_{\gamma,a}*_A\varphi\nu,\Psi\rangle_A=\langle P_{\gamma,a}*_A\varphi\nu,G_{a,2}*_ah_t\rangle_A $$ $$=\int_0^\infty\int_0^t\langle P_{\gamma,a,\tau} *_a(\varphi\nu)_{t-\tau},G_{a,2}*_a h_t\rangle_a d\tau dt$$ $$=\int_0^\infty\int_0^t\langle P_{\gamma,a,\tau} *_a G_{a,2}*_a(\varphi\nu)_{t-\tau},h_t\rangle_a d\tau dt=\langle P_{\gamma,a}*_A(G_{a,2}*_a(\varphi\nu)),h\rangle_A.$$
On the other hand, obviously if $C_{a,\gamma,p'}^{\frac{1}{p'}}(K)>0$ then there is a distribution $\nu$ supported on $K$ such that $\|P_{\gamma,a}*_A\nu\|_{p,A}<\infty$, that is $K$ is not removable.

\medskip

In Laplace case Proposition \ref{lp} can be extended to any $0< \gamma<2$. Indeed, all estimations hold in $0< \gamma<2$ case, except estimation of $II$, i.e. the terms which contain first space derivatives of $P_{\gamma,a}$. If $a=0$, the result immediately can be extended to $\gamma=1$, since in this case we have a singular integral, which can be estimated as in \cite[Lemma 3.5]{mp}. This method does not work when $0< \gamma<1$. The difficulty is caused by the fact that Bessel translation does not commute with product, i.e. $T^y(fg)(x)\neq T^y(f)(x)T^y(g)(x)$, cf. \eqref{tra1}. In the standard case there is no such problem.\\ On the other hand, although taking the Riesz kernel $I_\alpha(x):=\frac{c(\alpha)}{|x|^{n-\alpha}}$ instead of the Bessel one, we have that $-\Delta(I_\alpha*f)=I_\alpha*(-\Delta f)=I_{\alpha-2}*f$ (see \cite[page 118]{st}), but the estimation of the appropriate norms does not hold, i.e. $\|I_\alpha*f\|_p\not\le \|f\|_p$ (see \cite[Theorem 3.1.4]{ah}), so a lemma similar to Lemma \ref{l1} does not expected. Thus in the standard case we have the following theorem.

\medskip

\begin{theorem}\label{t2} Let $0< \gamma<2$ and $1< p <\infty $. Then for $a=0$, if $Z_{0,\gamma,p'}^{\frac{1}{p'}}(K)=0$, then $K\subset \mathbb{R}^{n+1}_+$ is removable for $L_{\gamma,0}$, and if $K\subset \mathbb{R}^{n+1}_+$ is removable for $L_{\gamma,0}$, then $C_{0,\gamma,p'}^{\frac{1}{p'}}(K)=0$.
\end{theorem}

\medskip

\noindent To prove Theorem \ref{t2}, we modify Lemma \ref{gammadelta}.

\medskip

\begin{lemma}\label{mgammadelta}
Let $\gamma \in (0,2)$, $1\le p<\infty$. Let $f\in \mathcal{S}$ and $\varphi \in C_0^\infty$. Then
\begin{equation}\label{1D}\left\|(-\Delta)^{\frac{\gamma}{2}}\varphi f\right\|_{p}\le c(n, \gamma)\left(\|\varphi\|_\infty  \|f\|_{W^{1,p}_{\Delta}}+\|\Delta\varphi\|_\infty \|f\|_p\right).\end{equation}
\end{lemma}

\proof Considering the integrals below in principal value sense, we have
$$-(-\Delta)^{\frac{\gamma}{2}}\varphi f=\int_{\mathbb{R}^n}\frac{(f\varphi)(x-y)-(f\varphi)(x)}{|y|^{n+\gamma}}dy$$ $$=\int_{\mathbb{R}^n}\varphi(x-y)\frac{f(x-y)-f(x)}{|y|^{n+\gamma}}dy+\int_{\mathbb{R}^n}f(x)\frac{\varphi(x-y)-\varphi(x)}{|y|^{n+\gamma}}dy=I+II.$$
Decomposing $I$ and $II$ to $|y|<1$ and $|y|\ge 1$ parts, as above, and denoting by $I=A+B$, $II=C+D$, repeating of the calculations above, we have
$$\|A\|_p\le \|\varphi\|_\infty |S^n_1|\int_0^1\frac{1}{r^{\gamma+1}}\|\Delta_{r,s,0}f(x)\|_pdr\le c(n,\gamma) \|\varphi\|_\infty \|\Delta f\|_p,$$
where Lemma \ref{gamma1} is applied.\\
Similarly
$$\|C\|_p\le \int_0^1\frac{1}{r^{\gamma+1}}\left(\int_{\mathbb{R}^n}\left|f(x)\int_{S^n_1}\varphi(x-r\Theta)-\varphi(x)dS(\Theta)\right|^pdx\right)^{\frac{1}{p}}dr$$ $$\le \int_0^1\frac{1}{r^{\gamma+1}}\|f\|_p|S^n_1|\|\Delta_{r,s,0}\varphi(x)\|_\infty dr\le  c(n,\gamma) \|\Delta \varphi\|_\infty \| f\|_p.$$
Of course,
$$\|B\|_p+\|D\|_p \le c(n,\gamma)\|\varphi\|_\infty \| f\|_p.$$
Thus triangle inequality ensures the result.

\proof (of Theorem \ref{t2})\\
The beginning of the proof up to  \eqref{eee} is the same as of Theorem \ref{lp}. After this, in view of Lemma \ref{mgammadelta}, we continue as
$$\|L_{\gamma,0}(\varphi(P_{\gamma,0}*\Psi)\|_{p'}$$ $$\le c(n, \gamma)\left(\|\varphi\|_\infty  (\|P_{\gamma,0}*\Psi\|_{p'}+\|\Delta(P_{\gamma,0}*\Psi)\|_{p'})+\|\Delta\varphi\|_\infty \|P_{\gamma,0}*\Psi\|_{p'}\right)$$ $$+\|\partial_t(\varphi (P_{\gamma,0}*\Psi)\|_p.$$
That is the terms which contain $\partial_{x_i}P_{\gamma,0}$ are missing. The remainder terms are already estimated, thus the proof of this direction is finished. The opposite direction is obvious as above.

\medskip

\vspace{5mm}

{\small{\noindent Both authors:\\
Department of Analysis and Operations Research,\\ Institute of Mathematics,\newline
Budapest University of Technology and Economics \newline
 M\H uegyetem rkp. 3., H-1111 Budapest, Hungary.

 \vspace{3mm}

\noindent \'A. P. Horv\'ath: g.horvath.agota@renyi.hu\\
M. Chegaar: mouna.chegaar@edu.bme.hu}

\end{document}